\documentclass[12pt]{amsart}
\usepackage{amssymb}
\usepackage{latexsym}
\usepackage{xy}
\input xypic
\newtheorem{theorem}{Theorem}%[section]
\newtheorem{proposition}[theorem]{Proposition}

\newtheorem{lemma}[theorem]{Lemma}

\def\Proof{\medskip\noindent{\bf Proof: }}

\def\Z{\mathbb{Z}}

\def\C{\mathbb{C}}

\def\R{\mathbb{R}}
\def\C{\mathbb{C}}

\def\F{\mathbb{F}}

\def\qed{\hfill$\square$\medskip}

\def\Zpk{\mathbb{Z}/p^{k}}
\def\Zpk1{\mathbb{Z}/p^{k-1}}

\newcommand{\rref}[1]{(\ref{#1})}

\newcommand{\beg}[2]{\begin{equation}\label{#1}#2\end{equation}}
\def\r{\rightarrow}

\def\sl2{\widetilde{SL_{2}(\Z)}}

\author{Daniel Kneezel and Igor Kriz}\thanks{Kriz is supported by NSA grant 08-1477.}
\title[Completing Verlinde Algebras]{Completing Verlinde Algebras}

\begin{document}

\maketitle

\begin{abstract}
We compute the completion of the Verlinde algebra of a simply connected
simple compact Lie group $G$ at the augmentation ideal of the representation ring. 
By results of Freed, Hopkins, Teleman and C.Dwyer and Lahtinen,
this gives a computation of (non-equivariant) twisted $K$-theory of the free loop space of $BG$.
\end{abstract}

\section{Introduction}
\label{s1}

The purpose of this note is to compute the completion of the Verlinde
algebra of a simply connected simple compact Lie group $G$ at the augmentation
ideal of the representation ring $R(G)$. In the special case of
of the symplectic groups, this was previously done in \cite{kw}. Denote the weight lattice of $G$ by
$\Pi^{*}=Hom(T,S^1)$ where $T$ is the maximal torus.
There is a unique inner product on $\Pi^*$ invariant under the action
of the Weyl group such that all long roots $u$ satisfy
$$\langle u,u\rangle=2.$$
The fundamental Weyl chamber $V$ is the set of all points in $\Pi^*\otimes\R$ with  
which all the positive roots have non-negative inner product. The {\em level $m$
fundamental alcove $A$} is the subset of $V$ of all points $x$
satisfying the relation 
\beg{ei1}{\langle u,x\rangle\leq m
}
where $a$ is the highest root in $V$ (cf. \cite{ps}). As usual, for weights $a,b$ we write 
$a\leq b$ when $b-a$ is non-negative on $V$.  A {\em level $m$ regular weight}
is a weight contained in the interior of $A$. Now the {\em level $m$ Verlinde algebra}
$V_m(G)$ is the quotient of $R(G)$ by the ideal $J_m$ of all elements $x$ such that
for every level $m$ regular weight $a$, 
\beg{ei2}{\phi_a(x)=0}
where $\phi_a$ is the composition of the inclusion
\beg{ei3}{\iota:R(G)\r \Z[\Pi^*]
}
with the map 
\beg{ei4}{\psi_a:\Z[\Pi^*]\r \C}
given by sending
\beg{ei5}{w\mapsto e^{2\pi i\langle w,a\rangle/m}.}
In this paper, we compute (as an abelian group) the completion of the Verlinde
algebra at the augmentation ideal $I$ of the ring $R(G)$. 

\vspace{3mm}
The definition of the Verlinde algebra, and our result, might at first appear artificial,
but our purely algebraic computation actually solves a problem in topology.
This is because of a theorem of Freed, Hopkins and Teleman \cite{fht},
which identifies the Verlinde algebra $V_m(G)$ as 
\beg{ei6}{K^{*}_{G,\tau}(G),}
the level $m$ equivariant twisted $K$-theory of $G$ with the conjugation action
of $G$. (The $*$ signifies that the non-trivial 
$K$-cohomology group is in the dimension
of parity equal to $dim(G)$; the twisted $K$-group in the dimension of the
other parity is $0$.)

\vspace{3mm}
Now it is well known (cf. \cite{kw}) that we have an equivalence
\beg{ei7}{G\times_{G} EG \simeq LBG
}
where on the left hand side of \rref{ei7}, the action of $G$ on $G$ is
by conjugation, and $\times_G$ means the quotient of the product by the
diagonal action of $G$, and the right hand side denotes the free loop space
on the bar construction of $G$.
Combining this with the completion theorem for twisted $K$-theory by
Dwyer and Lahtinen \cite{compl,compl1},
we see that the object we compute is actually the level $m$ (non-equivariant)
twisted $K$-theory of $LBG$. It is certainly nice to have an explicit computation
of these groups. As far as we know, this approach is the only known way of
computing them. A purely non-equivariant computation is not known (see
\cite{kw} for some comments on the difficulties of that approach). Therefore, the algebraic calculation
performed in the present paper is the key step to getting an answer,
which gives an instant source of examples of complete non-trivial computations
of non-equivariant twisted $K$-groups. 

\vspace{3mm}
The present note is organized as follows. In Section \ref{sr}, we state our
main results. In Section \ref{sl}, we will prove Theorem \ref{tvc}.
In Section \ref{slie}, we will discuss applications of Theorem \ref{tvc}
to the individual Lie group types, and prove the remainder of the
results of Section \ref{sr}.

\section{The main results}

\label{sr}

\begin{theorem}
\label{tvc}
The completion $(V_m(G))^{\wedge}_{I}$ is isomorphic, as an
abelian group, to the sum over primes $p\in\Z$ of $N(G,m,p)$
copies of $\Z_p$ where $N(G,m,p)$ is the number of
regular weights $a$ such that for every weight $w$, the denominator of
the rational number
$\langle w,a\rangle/m$ is a power of $p$.
\end{theorem}
\qed

The values of the numbers $N(G,m,p)$ are as follows. 
Let 
\beg{econv}{\parbox{3.5in}{$m=p^im^\prime$,
$n+1=p^\ell (n+1)^\prime$
$p$ does not divide $m^\prime$, $(n+1)^\prime$.}}

\vspace{3mm}
\noindent
{\bf Type $A$:} The number $N(A_n,m,p)$ is the number of tuples 
$$(b_1,...,b_n)\in\Z^n$$
such that
\beg{eai}{p^i>b_1>...>b_n>0,}
\beg{eaii}{(n+1)^\prime | (b_1+...+b_n).}

\vspace{3mm}
\noindent
{\bf Type $B$:} Assume $n>1$. For $p=2$, the number $N(B_n,m,2)$ is the number of tuples
\beg{ebtup1}{(b_1,...,b_n)\in \Z^n\cup (\Z+\frac{1}{2})^n}
such that
\beg{ebi}{b_1>...>b_n>0,}
\beg{ebii}{p^i>(b_1+b_2).}
For $p>2$, the number $N(B_n,m,p)$ is the number of tuples
\beg{ebtup2}{(b_1,...,b_n)\in\Z^n}
which satisfy \rref{ebi}, \rref{ebii} and 
\beg{ebiii}{2|(b_1+...+b_n).}

\vspace{3mm}
\noindent
{\bf Type $C$:} Assume $n>1$. For $p=2$, the number $N(C_n,m,2)$ is the number of tuples
\beg{ectup}{(b_1,...,b_n)\in\Z^n}
such that 
\beg{eci}{p^i>b_1>...>b_n>0.}
Explicitly,
$$N(C_n,m,2)= {2^i-1 \choose n}.$$
For $p>2$, the number $N(C_n,m,p)$ is the number of tuples \rref{ectup} which
satisfy \rref{eci} and
\beg{ecii}{2|b_i.}
Explicitly,
$$N(C_n,m,p)= {\frac{p^i-1}{2} \choose n}.$$

\vspace{3mm}
\noindent
{\bf Type $D$:} Assume $n>2$. For $p=2$, the number $N(D_n,m,2)$ is the number of tuples
\rref{ebtup1} such that
\beg{edi}{b_1>...>b_{n-1}>|b_n|,
}
\beg{edii}{p^i>b_1+b_2.}
For $p>2$, $N(D_n,m,p)$ is the number of tuples \rref{ebtup2}
which satisfy \rref{edi}, \rref{edii} and \rref{ebiii}.

\vspace{3mm}
\noindent
{\bf Type $G$:} For $p=3$, the number $N(G_2,m,3)$ is the number of tuples
\beg{egtup}{(b_1,b_2)\in \Z^2}
which satisfy
\beg{egi}{2b_2>b_1>b_2>0,
}
\beg{egii}{p^i>b_1.}
For $p\neq 3$, $N(G_2,m,p)$ is the number of tuples \rref{egtup} which satisfy
\rref{egi}, \rref{egii} and 
\beg{egiii}{3|(b_1+b_2).}

\vspace{3mm}
\noindent
{\bf Type $F$:} For $p=2$, the number $N(F_4,m,p)$ is the number of
tuples \rref{ebtup1} for $n=4$ such that
\beg{efi}{b_2>b_3>b_4>0,}
\beg{efii}{b_1>b_2+b_3+b_4,}
\beg{efiii}{p^i>b_1+b_2.}
For $p>2$, $N(F_4,m,p)$ is the number of tuples \rref{ebtup2} for $n=4$
which satisfy \rref{efi}, \rref{efii}, \rref{efiii} and \rref{ebiii} for $n=4$.

\vspace{3mm}
\noindent
{\bf Type $E_8$:} The number $N(E_8,m,p)$ is the number of tuples
\beg{etupe8}{(b_1,...,b_8)\in \Z^8\cup(\Z+\frac{1}{2})^8
}
such that
\beg{ee8i}{b_2>...>b_7>|b_8|,}
\beg{ee8ii}{b_1>b_2+...+b_7-b_8,}
\beg{ee8iii}{p^i>b_1+b_2,}
and \rref{ebiii} for $n=8$.

\vspace{3mm}
\noindent
{\bf Type $E_7$:} For $p=2$, the number $N(E_7,m,2)$ is the number of tuples
\beg{etupe7}{(b_1,...,b_7)\in \frac{1}{\sqrt{2}}\Z\times(\Z^6\cup(\Z+\frac{1}{2})^6)
}
such that
\beg{ee7i}{b_2>...>b_6>|b_7|,}
\beg{ee7ii}{b_1\sqrt{2}>b_2+...+b_6-b_7,}
\beg{ee7iii}{p^i>\sqrt{2}b_1,}
\beg{ee7iv}{2|(\sqrt{2}b_1+b_2+...+b_6-b_7).}
For $p>2$, the number $N(E_7,m,2)$ is the number of tuples 
\rref{etupe7} where 
\beg{etupe72}{2b_i\equiv \sqrt{2}b_1\mod 2\; \text{for $i=2,...,7$}
}
such that \rref{ee7i}, \rref{ee7ii}, \rref{ee7iii} and
\beg{ee7v}{2|(b_2+...+b_7).}

\vspace{3mm}
\noindent
{\bf Type $E_6$:} For $p=3$, the number $N(E_6,m,3)$ is the number of tuples 
\beg{etupe63}{(b_1,...,b_6)\in (\frac{1}{\sqrt{3}}\Z\times\Z^5)\cup 
(\frac{1}{\sqrt{3}}(\Z+\frac{1}{2})
\times(\Z+\frac{1}{2})^5)
}
such that
\beg{ee6i}{b_2>...>b_5>|b_6|,}
\beg{ee6ii}{\sqrt{3}b_1>b_2+...+b_5-b_6,}
\beg{ee6iii}{p^i>(\sqrt{3}b_1+b_2+...+b_6)/2,}
\beg{ee6iv}{2|(\sqrt{3}b_1+b_2+...+b_6).}
For $p\neq 3$, the number $N(E_6,m,2)$ is the number of tuples
\beg{etupe6}{(b_1,...,b_6)\in (\sqrt{3}\Z\times\Z^5)\cup (\sqrt{3}(\Z+\frac{1}{2})
\times(\Z+\frac{1}{2})^5)
}
such that \rref{ee6i}, \rref{ee6ii}, \rref{ee6iii}, \rref{ee6iv}.

\vspace{3mm}

{\bf Remark:} As already remarked
in \cite{kw} in the case of $Sp(n)$, it
is interesting to note that completion 
at the augmentation ideal of the Verlinde algebra does 
not preserve the level-rank
duality \cite{bowr}. This indicates that twisted $K$-theory of $LBG$ is perhaps
not as natural an object in mathematical physics as the Verlinde algebra
(which can be interpreted as the fusion ring of certain conformal field
theories known as WZW models).

\section{The fundamental lemmas}
\label{sl}

\begin{lemma}
\label{laug}
There exists a number $N$ such that $N.1\in J_m$ where $J_m$
is in the augmentation
ideal of $V_m(G)$ (i.e. the image $\overline{I}$ of the ideal $I$ in $V_m(G)$).
\end{lemma}

\Proof
If an element $w\in J_m$ has non-zero augmentation $N$, then clearly
$N\in \overline{I}$. Thus, it suffices to show that
\beg{els1}{J_m\nsubseteq I.
}
Suppose thus that \rref{els1} is false, i.e. that $J_m\subseteq I$. Then 
the augmentation $\epsilon:R(G)\r \Z$ must factor through $V_m(G)$.
After tensoring with $\C$, the product of the maps $\phi_a$ over
regular weights $a$ induces an iso of $V_m(G)\otimes \C$ with
the product of rings
\beg{els2}{\prod_a\C} 
over regular weights $a$. Thus, the complexified augmentation must factor
through a map of rings from \rref{els2} to $\C$. But one easily sees that
the set of maps of rings from \rref{els2} to $\C$ is the set of projections
(since such maps must preserve idempotents). So, such a factorization 
would mean that the augmentation is the map $\phi_a$
for some $a$. However,
averaging over complex conjugates of roots of unity, we then get that
the augmentation must coincide with the map which is $1$ on
any weight $w$ satisfying $m|\langle w,a\rangle$, and $0$ on other weights.
This is clearly false, as any weight can occur in a representation.
\qed

By Lemma \ref{laug}, the completion of the Verlinde algebra is
a finitely generated module over the ring of $N$-adic numbers $\Z_N$,
and hence a finite sum of finitely generated $\Z_p$-modules
for finitely many primes $p$. Our main tool is this: the
canonical inclusion from $V_m(G)$ into the ring \rref{els2}. Then this
map factors through a finite integral extension $\Z^\prime$ of $\Z$
(just adjoin all the necessary roots of unity). Then, we obtain a map
of $\Z^\prime$-algebras
\beg{esf1}{V_m(G)\otimes \Z^\prime \r \prod_a\Z^{\prime}
}
where the product on the right hand side is over level $m$ regular weights.
Denote the left hand side of \rref{esf1} by $R$, and the right hand side by $R^\prime$.
Thus, we have a short exact sequence of $R$-modules
\beg{esf2}{0\r R\r R^\prime\r F\r 0
}
where $F$ is finite (since it is finitely generated and vanishes after
tensoring with $\C$). Our strategy is to compute the completion of $R^\prime$
instead of $R$, with the idea that the finite cokernel \rref{esf1} will not
make a difference in the completion. This is justified by the following
result.

\vspace{3mm}
\begin{theorem}
\label{t1}
The sequence
\beg{etf1}{0\r \varprojlim R/I^n \r \varprojlim R^\prime/I^n \r \varprojlim F/I^n\r 0
}
induced by \rref{esf2} is exact. Further, the first non-zero term of \rref{etf1} 
is torsion free.
\end{theorem}

\Proof
First note that there exists a number $N$ such that $NR^\prime\subseteq R$.
We conclude that 
$$[I^nR^\prime :I^nR]\leq [I^nR^\prime:NI^nR^\prime]\leq N^{k}$$
where $k$ is the number of cyclic summands of $R^\prime$. 
Thus, we obtain a short exact sequence
\beg{etf2}{0\r Q_n\r R/I^nR\r R/((I^nR^\prime)\cap R)\r 0}
where $Q_n$ are finite of bounded order. It follows that the inverse limit
of \rref{etf2} is a short exact sequence of the form
\beg{etf3}{0\r Q\r \varprojlim R/I^nR\r \varprojlim R/((I^nR^\prime)\cap R)\r 0
}
where $Q$ is finite.  If we know that the middle term \rref{etf3}
is torsion free, then $Q=0$. Now by \cite{am}, Proposition 10.3, we have a short
exact sequence
\beg{etf4}{0\r \varprojlim R/((I^nR^\prime)\cap R)\r
\varprojlim R^\prime/I^nR^\prime \r \varprojlim T/I^n F.
}
This, together with \rref{etf3} then implies the statement.

To prove that $\varprojlim R/I^nR$ is torsion free, let $p_i$ be the composition
of \rref{esf1} with the projection of the right hand side to the first $i$ factors.
Then $p_i$ is a map of rings, so its kernel $K_i$ is an ideal in $R$, and
we have a finite filtration of $R$ by ideals
$$R\supset K_1\supset K_2\supset...\supset 0$$
where $K_i/K_{i+1}$ is isomorphic to an ideal of $\Z^\prime$.
By Proposition 10.3 of \cite{am}, the completion of $R$ with respect
to the augmentation ideal $I$ is filtered by the completions of $K_i$
with respect to the filtrations $F(i)_n=K_i\cap I^n$, where
the associated graded object is the sum of completions of $K_i/K_{i+1}$
by the image of the filtration $F(i)$.

It suffices to prove then that the completion of $K_i/K_{i+1}$ by
the image of the filtration $F(i)$ is non-torsion. But since
$$I^n\cdot K_i\subseteq F(i)_n,$$
there is an onto map
\beg{etff2}{(K_i/K_{i+1})^{\wedge}_{I}\r (K_i/K_{i+1})^{\wedge}_{F(i)}.
}
Therefore, the right hand side of \rref{etff2} can have non-zero torsion
only when it is non-zero, which is when the left hand side of \rref{etff2}
is non-zero.  But $\Z^\prime$ is a Dedekind domain, so a completion of
any of its ideal by another non-zero ideal $I$ is isomorphic to the
completion of $\Z^\prime$ by $I$. (The image of the ideal $I$ is
non-zero by Lemma \ref{laug}.)

In other words, taking the direct 
sum over $i$, the left hand side of \rref{etff2}
becomes the completion of the right hand side of \rref{esf1}, 
the right hand side of \rref{etff2} is the associated graded of the completion
of the left hand side of \rref{esf1}. Therefore, if one of the right hand sides
of \rref{etff2} contains torsion at a prime $\pi$ in $\Z^\prime$, then
the $\Z^{\prime}_{\pi}$-rank of the completion of the right hand side
of \rref{esf1} is greater than the $\Z^{\prime}_{\pi}$-rank of the completion
of the left hand side of \rref{esf1}, which contradicts \rref{etf3} 
and \rref{etf4}.
\qed

\vspace{3mm}
Thus, it suffices to compute the completion of the $R$-module $R^\prime$
at the image of the ideal $I$. Since $R^\prime$ is a direct sum of $R$-modules,
it suffices to compute the completion of each summand separately. We will
denote the $a$-th summand by $\Z^{\prime}_{a}$. 

\begin{proposition}
\label{p1}
$(\Z^{\prime}_{a})^{\wedge}_{I}$ is the sum of $\Z^{\prime}_{p}$ over all
primes $p\in \Z$ such that the denominator of the rational
number $\langle w,a\rangle/m$
is a power of $p$.
\end{proposition}

This immediately implies Theorem \ref{tvc}.
To prove Proposition \ref{p1}, we will need the following lemma:

\begin{lemma}
\label{lp11}
Let $\pi$ be a prime of $\Z^{\prime}$ which lies above a prime
$p\in\Z$. Then a primitive $m$'th root of
unity $\zeta_m\in \Z^\prime$ satisfies 
\beg{ezeta}{\pi|(\zeta_m-1)}
if and only if $m$ is a power of $p$.
\end{lemma}

\Proof
If $m=p^j$, then $x=\zeta_m$ satisfies $x^{p^j}=1$ in the residue field,
which clearly implies $x=1$. Conversely, when $m$ is not a power of $p$,
we may as well assume that $m$ is relatively prime to $p$ (since \rref{ezeta}
implies the same statement with $m$ replaced by any of its factors). 
But when $m$ is relatively prime to $p$, then $\zeta_m$ is a root
of the polynomial 
\beg{ezeta1}{x^{m-1}+x^{m-2}+...+1,} 
and hence the same is true
in the residue field. But the polynomial \rref{ezeta1} does not have
the root $1$ in $\overline{\F_p}$.
\qed

\vspace{3mm}
\noindent
{\bf Proof of Proposition \ref{p1}:}
It suffices to compute the completion of $\Z^{\prime}_{a}$ at
every ideal $(I,\pi)$ where $\pi$ is a prime of $\Z^{\prime}$.
This is $\Z_{\pi}$ when every element of $I$ is divisible
by $\pi$, and $0$ else. The condition for the completion
being nonzero then is
\beg{enon1}{\pi|(u-dim(u))
}
for the image $u$ of every element of $R(G)$ in $\Z_a$.
We claim that this is equivalent to \rref{enon1} holding for
every weight $u$! In effect, $\Z[\Pi^*]$ is an
integral extension of $R(G)=\Z[\Pi^*]^{W(G)}$,
so every weight $x=u$ is a root of a polynomial
\beg{enon2}{x^n+w_1x^{n-1}+....+w_n
}
where $w_i\in R(G)$ (the other roots are $W(G)$-conjugate weights). 
Let $dim(w_i)=m_i$. When \rref{enon1}
holds for any weight, clearly it holds in particular for
any element of the representation ring. When \rref{enon1}
holds for any element of the representation ring,
the polynomial \rref{enon2} is equal to
\beg{enon3}{x^n +m_1 x^{n-1}+...+m_n 
}
over the residue field. But then all the roots of \rref{enon3} are of the form
$dim(u^\prime)=1$ where $u^\prime$ are weights
conjugate to $u$, since the augmentation is a ring homomorphism. 
Thus, \rref{enon3} is in fact equal to $(x-1)^n$ over the residue field,
which shows that the image of the weight $u$ (which was a root of
\rref{enon2}) in the residue field must be equal to $1$.

Thus, we know that the completion of $\Z_a$ is nonzero
if and only if
\beg{enon4}{\pi|(u-1)
} 
for every weight $u$. Now $u$ is the $\langle u,a\rangle$'th power of
$\zeta_m$, the $m$'th primitive root of unity. Thus, the statement follows from 
Lemma \ref{lp11}.
\qed

\vspace{3mm}
\section{Discussion of the individual Lie group types}
\label{slie}

The purpose of this section is to deduce from Theorem \ref{tvc}
the computations of the numbers $N(G,m,p)$ in the individual types
of simple simply connected compact Lie groups. Recall convention 
\rref{econv}. For basic facts about compact Lie groups and
their representation, see \cite{fh}.

\vspace{3mm}
\noindent
{\bf Type $A$:} The weight lattice $\Pi^*$ for $A_n$ ($SU(n+1)$) can be identified
with the quotient of $\Z^{n+1}$ by the subgroup $\langle 1,...,1\rangle$.
The inner product is such that the orthogonal projection to 
\beg{earoot}{\{(x_1,...,x_{n+1})\in \R^{n+1}|\sum x_i=0\}}
is an isometry (where we consider the induced inner product from $\R^{n+1}$ on
\rref{earoot}). The roots are vectors of the form
$$(0,...,0,1,0,...,0,-1,0,....,0)$$
(the numbers $1$ and $-1$ are in two arbitrary different places). Thus,
a fundamental Weyl chamber can be selected as the set of all weights
\beg{eaft}{a=(a_1,....,a_n,0)} 
where
\beg{eafund}{a_1\geq....\geq a_n\geq 0.
}
The highest root is 
$$(1,0,...,0,-1),$$
so the regular weights are those tuples \rref{eaft} where $a_i\in \Z$, and
\beg{aaalc}{m>a_1>...>a_n>0.
}
The weight lattice is generated by the weights
$$(0,...,0,1,0,...,0)=(-1/(1+n),....,-1/(1+n),n/(1+n),-1/(1+n),....),$$
so the condition of Theorem \ref{tvc} is that the denominator of
each of the numbers
\beg{easum1}{\sum_{i=1}^{n} -a_i/(1+n),}
\beg{easum2}{\sum_{i=1}^{n} -a_i/(1+n)+a_j}
for any fixed $j$ is a power of $p$. Clearly, we may
replace the numbers \rref{easum2} with the $a_i$'s.
Now put
\beg{eab}{b_i:=a_i/m^\prime.}
Translating our conditions fro the $a_i$'s to the $b_i$'s gives the condition
stated in Section \ref{sr}.

\vspace{3mm}
\noindent
{\bf Type $B$:} The weight lattice of $B_n$ ($Spin(2n+1)$ - we assume here $n>1$,
for $n=1$, the group is $SU(2)$) is the set of all tuples
\beg{ebw}{\text{$(a_1,...,a_n)$
where all $a_i$'s are either integers or integers plus $1/2$.}} 
The inner product
is induced from $\R^n$. The long roots are all vectors of the form
$$(0,...,0,\pm 1,0,....,0,\pm 1,0,...,)$$
(the two nonzero entries are in arbitrary different places and the
signs are not required to match), and the short roots are of the
form
$$(0,...,0,\pm 1,0,...,0).$$
Thus, the fundamental Weyl chamber can be chosen as the set of all weights
\rref{ebw} such that \rref{eafund} holds. The highest root is then 
$$(1,1,....,0,0),$$
so the regular weights are all those satisfying
\beg{ebrw}{\text{$a_1>...>a_n>0$ and $m>a_1+a_2$.}}
The generating weights can be taken as
\beg{ebw1}{(0,...,0,1,0,...,0),}
\beg{ebw2}{(1/2,...,1/2).}
Putting again \rref{eab} gives the conditions of Section \ref{sr}. Note
that when $p=2$, the denominator $2$ does not matter in the condition,
so the weight \rref{ebw2} does not contribute to the condition, and also 
\rref{ebw} remains the same with $a_i$'s replaced by $b_i$'s. On
the other hand, for $p>2$, a denominator $2$ violates the condition
of Theorem \ref{tvc}, which is why the numbers $b_i$ must be integers,
and their sum must be even to get an integral inner product with the
weight \rref{ebw2}.

\vspace{3mm}
\noindent
{\bf Type $C$:} The weight lattice of $C_n$ (usually denoted $Sp(n)$ although
sometimes also $Sp(2n)$) is the set of all tuples 
\beg{ecw}{\protect\parbox{3.5in}{$(a_1,...,a_n)$ such that $\protect
a_i\sqrt{2}\in \Z$.}}
The inner product is again induced from $\R^n$. The long roots are
$$(0,...,0,\sqrt{2},0,...,0),$$
the short roots are
$$(0,...,0,\pm 1/\sqrt{2},0,....,0,\pm 1/\sqrt{2},0,...,)$$
(the nonzero entries are in arbitrary places and the signs
are independent). Thus, the fundamental chamber can be selected
as the set of all weights \rref{ecw} where \rref{eafund} holds.
The highest root is
$$(\sqrt{2},0,....,0),$$
so the regular weights are those satisfying
\beg{ecreg}{m/\sqrt{2}>a_1>...>a_n>0.
}
The weights are generated by
$$(0,...,0,1/\sqrt{2},0,...,0).$$
To translate the condition of Theorem \ref{tvc}, we must now put
\beg{ecb}{b_i:=a_i\sqrt{2}/m^\prime.
}
Clearly, we obtain the conditions listed in Section \ref{sr}. Note in 
particular again that when the $b_i$'s are integers, inner product
with a generator weight can produce a denominator of $2$, which
can be neglected at $p=2$, but violates the condition for $p>2$. This
is why for $p>2$, the numbers $b_i$ are required to be even. 

\vspace{3mm}
\noindent
{\bf Type $D$:} For the group $D_n$ ($Spin(2n)$),
the weight lattice is the same as for $B_n$ (with the same
inner product), but only the
long roots of $B_n$ are roots of $D_n$ (we assume here $n>2$ - for $n=2$ 
the group is $SU(2)\times SU(2)$). Thus, the fundamental
Weyl chamber is the set of all weights \rref{ebw} 
such that
\beg{edch}{a_1\geq...\geq a_{n-1}\geq|a_n|.
}
The highest root remains the same as for $B_n$, so 
the regular weights are those satisfying
\beg{edreg}{a_1>...>a_{n-1}>|a_n|\; m>a_1+a_2.
}
Since the weights are the same as for 
$B_n$, clearly, using \rref{eab}, the divisibility conditions for the
$b_i$'s are the same as in the $B_n$ case, while \rref{edreg}
translates to \rref{edi}, \rref{edii}.

\vspace{3mm}
\noindent
{\bf Type $G$:} The weight lattice of $G_2$ is the same as for $A_2$.
The long roots are the same as for $A_2$, but there are also short
roots
$$(\pm 1,0,0), (0,\pm1,0),(0,0,\pm 1).$$
Thus, the fundamental Weyl chamber can be
chosen as the set of all weights \rref{eaft}
such that 
$$2a_2\geq a_1\geq a_2\geq 0.$$
The highest root is the same as for $A_2$, so the fundamental chamber
is the set of all weights \rref{eaft} such that
$$2a_2>a_1>a_2>0,\; m>a_1.$$
Using again \rref{eab}, we obtain the conditions \rref{egi}, \rref{egii},
and the same divisibility conditions as for $A_2$.

\vspace{3mm}
\noindent
{\bf Type $F$:} The weight lattice of $F_4$ is the same as for $B_4$. The
roots are those of $B_4$, plus the roots
$$(\pm 1/2,\pm 1/2,\pm 1/2,\pm 1/2)$$
(the signs are independent). Thus, the fundamental Weyl chamber
can be selected as the set of all weights \rref{ebw} such that
$$\geq a_2\geq a_3\geq a_4\geq 0,\; a_1\geq a_2+a_3+a_4.$$
The highest root is the same as for $B_4$. Thus, the regular weights
are those weights \rref{ebw} where
\beg{efreg}{a_2>a_3>a_4>0,\; a_1>a_2+a_3+a_4,\; m>a_1+a_2.}
Using \rref{eab}, \rref{efreg} translates to \rref{efi}, \rref{efii}, \rref{efiii},
and the divisibility conditions remain the same as for $B_4$.

\vspace{3mm}
\noindent
{\bf Type $E_8$:} The weight lattice of $E_8$ consists of all $8$-tuples
\beg{ewe8}{\parbox{3.5in}{$(a_1,...,a_8)$ such that $a_i$ are all integers or
integers plus $1/2$, and their sum is even.}
}
(The evenness requirement comes from the fact that $E_8$ does not
actually contain $Spin(16)$, only $SO(16)$.) The roots are all the
elements of the form
$$(0,...,0,\pm 1,0,...,0,\pm 1,0,...,0)$$
where the nonzero entries are in arbitrary places and the signs are independent,
and elements of the form
$$(\pm 1/2,...,\pm 1/2)$$
where the signs are arbitrary such that the number of minus signs is even. 
Thus, a fundamental Weyl chamber can be chosen as the set of all
weights
\beg{eche8}{a_2\geq a_3\geq...\geq a_7\geq |a_8|,\; a_1\geq a_2+...+a_7-a_8
}
and the highest root is
$$(1,1,...,0,0),$$
so the regular weights are all weights satisfying
\beg{erege8}{a_2>...>a_7>|a_8|, \; a_1>a_2+...+a_7-a_8,\; m>a_1+a_2.
}
The weights are generated by
\beg{ew1e8}{(0,...,0,1,-1,0,...,0)
}
and
\beg{ew2e8}{(1/2,...,1/2).}
So, if we use \rref{eab}, we obtain the conditions \rref{ee8i}, \rref{ee8ii},
\rref{ee8iii}. Note that due to the unimodularity of the
$E_8$ lattice, the divisibility condition becomes vacuous. (Note
also that $E_8$ is the only case in which there are no
exceptional primes.)

\vspace{3mm}
\noindent
{\bf Type $E_7$:} The weight lattice of $E_7$ can be identified
with the orthogonal projection of the $E_8$ lattice \rref{ewe8} to the orthogonal
complement of a fixed root of $E_8$. Choosing the $E_8$  root
$$(1,1,0,...,0),$$
the weights are all elements 
\beg{ew1e7}{(a_1/\sqrt{2},-a_1/\sqrt{2},a_2,...,a_7)
}
where
\beg{ew2e7}{\parbox{3.5in}{$a_1\in \frac{1}{\sqrt{2}}\Z$, $a_2,...,a_7$ are
all either all integers or integers plus $1/2$ and $\sqrt{2}a_1+a_2+...a_6-a_7$ is even.}}
The roots are 
$$(0,...,0,\pm 1,0,...,0,\pm 1,0,...,0)$$
where the non-zero positions are in arbitrary places other than the first
coordinate, and the signs are independent, 
$$(\pm\sqrt{2},0,....,0),$$
and
$$(\pm 1/\sqrt{2},\pm 1/2,...,\pm 1/2)$$
where the signs are independent such that the number of minus signs is odd. 
Therefore, a fudamental Weyl chamber can be selected as the set of all
weights \rref{ew1e7}, \rref{ew2e7} such that
\beg{efunde7}{a_2\geq...\geq a_6\geq |a_7|,\; a_1/\sqrt{2}\geq a_2+....+a_6-a_7.
}
The highest root is 
$$(\sqrt{2},0,...,0),$$
so the regular weights are subject to the conditions
\beg{erege7}{a_2>...>a_6>|a_7|,\; m>a_1\sqrt{2}>a_2+...+a_6-a_7.
}
The weights are generated by 
\beg{ewg1e7}{(\sqrt{2},0,...,0),}
\beg{ewg2e7}{(1/\sqrt{2},0,...,0,1,0,...,0),}
\beg{ewg3e7}{(0,1/2,...,1/2).}
Now using \rref{eab} leads to the divisibility conditions in Section \ref{sr}.
Note that for $p=2$, a factor of $2$ in the denominator of $b_i, i>1$ 
or $b_1\sqrt{2}$ does not 
violate the condition of Theorem \ref{tvc}, so the condition of Theorem \ref{tvc}
with respect to the weights \rref{ewg1e7}, \rref{ewg2e7}, \rref{ewg3e7}
are automatically satisfied by the conditions \rref{etupe7}, \rref{ee7iv},
which in turn follow from \rref{ew1e7}, \rref{ew2e7}. For $p>2$, 
on the other hand, a $2$ in the denominator of $b_1/\sqrt{2}+b_i, i>1$ 
and $b_2+...+b_7$
violates the condition of Theorem \ref{tvc}, which forces \rref{etupe72}
and \rref{ee7v}. The condition of Theorem \ref{tvc}
with \rref{ewg1e7} is automatic.

\vspace{3mm}
\noindent
{\bf Type $E_6$:}
The weight lattice of $E_6$ can be identified with the orthogonal projection
of the $E_8$ lattice to the orthogonal complement of the $E_8$  roots
$$(1,-1,0,...,0), \; (0,1,-1,0,...,0).$$
This gives weights
\beg{ew1e6}{(a_1/\sqrt{3},a_{1}/\sqrt{3},a_1/\sqrt{3},a_2,...,a_6) 
}
where
\beg{ew2e6}{\parbox{3.5in}{$a_1\sqrt{3},a_2,...,a_6$ are all either integers or
integers plus $1/2$}
}
and
\beg{ew3e6}{2|(\sqrt{3}a_1+a_2+...+a_6).}
The roots are  
$$(0,...,0,\pm 1,0,...,0,\pm 1,0,...,0)$$
where the non-zero positions are in arbitrary places other than the first
coordinate, and the signs are independent, 
and
$$(\pm \sqrt{3}/2,\pm 1/2,...,\pm 1/2)$$
where the signs are independent such that the number of minus signs is even. 
A fundamental Weyl chamber can be selected as the set of all weights
\rref{ew1e6}, \rref{ew2e6}, \rref{ew3e6} satisfying
\beg{efunde6}{a_2\geq ...\geq a_5\geq |a_6|,\; \sqrt{3}a_1\geq a_2+...+a_5-a_6.
}
The highest root is 
$$(\sqrt{3}/2,1/2,...,1/2),$$
so the regular weights are those satisfying the condition
\beg{erege6}{a_2>...a_5>|a_6|,\; \sqrt{3}a_1>a_2+...+a_5-a_6,\; m>\sqrt{3}a_1+a_2+...a_6)/2.
}
The weights are generated by 
\beg{ewg1e6}{(2/\sqrt{3},0,...,0),}
\beg{ewg2e6}{(1/\sqrt{3},0,...,0,1,0,...,0),}
\beg{ewg3e6}{(\sqrt{3}/2,1/2,...,1/2).}
Using \rref{eab} gives the divisibilty conditions in Section \ref{sr}.
Note that at $p=3$, a factor of $3$ in the denominator does not matter,
so the defining conditions for weights \rref{ew1e6}, \rref{ew2e6}, \rref{ew3e6}
imply the divisibility conditions with respect to the weights
\rref{ewg1e6}, \rref{ewg2e6}, \rref{ewg3e6}. On the other hand, at $p\neq 3$,
a denominator of $3$ in $\sqrt{3}a_1$ immediately violates the condition
of Theorem \ref{tvc} with respect to the weight \rref{ewg1e6}, which forces
\rref{etupe6}. This together with the defining conditions \rref{ew2e6} and
\rref{ew3e6} forces the condition of Theorem \ref{tvc} for the weights
\rref{ewg2e6}, \rref{ewg3e6}.

\end{document}